\documentclass[reqno,12pt]{amsart}

\usepackage[margin=2.5cm]{geometry}
\usepackage[T1]{fontenc}
\usepackage[latin1]{inputenc}
\usepackage[english]{babel}
\usepackage{amsmath,amssymb,amsthm}
\usepackage{enumitem}
\usepackage{mymacros}
\usepackage{hyperref}
\usepackage{xcolor}
\usepackage{graphicx}

\DeclareMathOperator{\capa}{cap}

\renewcommand{\MR}[1]{}
\DeclareMathOperator{\hausdorff}{dim}
\title{Riesz capacities of a set due to Dobi\'nski}

\author{Nicola Arcozzi}
\address{Dipartimento di Matematica, Universit\`a di Bologna, 40126, Bologna, Italy}
\email{nicola.arcozzi@unibo.it}

\author{Nikolaos Chalmoukis}
\address{Dipartimento di Matematica, Universit\`a di Bologna, 40126, Bologna, Italy}
\email{nikolaos.chalmoukis2@unibo.it}
\thanks{This research was supported by the Ministry of Science and Higher Education of the Russian Federation, agreement No. 075-15-2019-1619}


\begin{document}
\begin{abstract}
     We study the Riesz $(a,p)$-capacity of the so called Dobi\'nski set. We characterize the values of the parameters $a$ and $p$ for which the $(a,p)$-Riesz capacity of the Dobi\'nski set is positive. In particular we show that the Dobi\'nski set has positive logarithmic capacity, thus answering a question of Dayan, Fernand\'ez and Gonz\'alez. We approach the problem by considering the dyadic analogues of the Riesz $(a,p)$-capacities which seem to be better adapted to the problem.
\end{abstract}
\maketitle
\section{Introduction and main results}

In a series of two papers \cite{Dobinski1876, Dobinski1877} Dobi\'nski claims that the following identity is true 
\[ \prod_{n\geq 0} (\tan 2^n \pi x )^{2^{-n}} = (2\sin \pi x)^2, \]
for all real numbers $x\in[0,1]$ which are not dyadic rationals. As it has been already noted in \cite{Agnew1947} and explained in detail in a recent paper \cite{Dayan2021} the situation is not quite so simple. In fact, if we consider the same identity with absolute values, 
\[ \prod_{n\geq 0} |\tan 2^n \pi x |^{2^{-n}} = (2\sin \pi x)^2,  \]
so as to avoid issues of defining the powers of negative numbers, in \cite{Agnew1947} the authors prove that the identity holds if and only if $x$ does not belong to the so called {\it Dobi\'nski set} $\cD$.  To define $\cD$ let $x\in[0,1]$ be a real number with dyadic expansion $x=(0.a_1a_2\dots)_2$ and for $n\geq 1$ let $s_n(x) = \max\{r\in \bN: a_n = a_{n+1}= \dots =a_{n+r} \}$ \footnote{Let $s_n(x)=\infty$ if $x$ is a dyadic rational}. Then,

\[ \cD:= \{ x\in [0,1]: \limsup_{n\to \infty} \frac{s_n(x)}{2^n} > 0  \}. \]

So the Dobi\'nski set  comprises real numbers which can be approximated ``exceedingly well'' by dyadic rationals on every scale. Related problems of diophantine approximation by dyadic rationals have been considered in \cite{Nilsson2009}.

In a recent work \cite{Dayan2021} Dayan, Fernand\'ez and Gonz\'alez prove, among other results, that $\cD$ has Hausdorff dimension $0$ and {\it logarithmic} Hausdorff dimension $1$. Their techniques are primarily based on the mass transference principle of Beresnevich and Velani \cite{Beresnevich2006} which allows one to transfer measure theoretic statements for $\limsup$ subsets of $\bR^n$ to statements about Hausdorff measure.  
The Hausdorff dimension is a precise way to talk about the size of a subset of $\bR^n$. Another way to measure the size of subsets of $\bR^n$ is by some kind of {\it capacity}. From now on it will be convenient to consider $\cD$ as a subset of the unit circle $\bT$ in $\bR^2$ via the usual correspondence $x \mapsto e^{2\pi i x} $. Let   $0<a<1$ and $f $  a positive measurable function on $\bT$ the {\it $a$-Riesz potential of $f$} is defined as 
\[ \cI_a f(x):= \int_{\bT}\frac{f(y)dy}{|x-y|^{1-a}}, \quad x\in \bT, \]
where $dy$ is the normalized Lebesgue measure on $\bT.$
Finally let $1<p<\infty$ and $a$ as before. The {\it $(a,p)$-Riesz capacity} of a Borel subset  $E$ of $\bT$ is defined as 
\[ R_{a,p}(E):= \inf \big\{\int_{\bT}f(y)^p dy :f \geq 0\,\, \text{and}\,\, \cI_a f \geq 1 \,\, \text{on} \,\, E \big\}. \]
We shall refer to the capacities $R_{a,p}$ as {\it linear} if $p=2$ and as {\it non-linear} if $p\neq 2$. We shall also assume that $ap\leq 1$, otherwise singletons have positive capacity.
There exists a remarkable relation between Hausdorff dimension, which we will denote by $\hausdorff$, and linear Riesz capacities for a Borel set $E\subseteq \bT$ established by Frostman \cite[Corollary 5.1.14]{Adams1996}
\begin{equation}\label{eq:hausdorff} \hausdorff E =\sup \{1-a p: R_{a,p}(E)>0 \}.\end{equation}
This fact, together with the standard comparison results between capacities \cite[Section 5.5]{Adams1996} and the fact that the Dobi\'nski set has vanishing Hausdorff dimension implies that $R_{a,p}(\cD)=0$ when $ap<1$. In the same work where the Hausdorff dimension of $\cD$ was studied the authors ask whether also $R_{\frac{1}{2},2}(\cD) = 0 $ or not \cite[Section 5]{Dayan2021}. In fact they formulate the question in terms of {\it logarithmic capacity } in the complex plane but it is well known that for subsets of $\bT$, logarithmic capacity is bounded below and above by the Riesz $(\frac{1}{2},2)$-capacity \cite[Corollary 2.6]{Chalmoukis2020}.

We have been able to answer the above question for all Riesz $(\frac{1}{p},p)$-capacities. 

\begin{thm}\label{thm:main}
Let $\cD$ be the Dobi\'nski set and $p>1$. Then, 
\[ R_{\frac{1}{p},p}(\cD) =  R_{\frac{1}{p},p}(\bT ) \,\,\, \text{if} \,\,\, 1<p\leq2 \,\,\, \text{and} \,\,\,  R_{\frac{1}{p},p}(\cD) = 0 \,\,\, \text{if} \,\,\, p > 2. \]
\end{thm}

Somewhat surprisingly the capacity of $\cD$ exhibits a jump from full to $0$ at the critical value $p=2$. It should be mentioned that this statement implies, via \eqref{eq:hausdorff}, that the Hausdorff dimension of $\cD$ is $0$  and that the logarithmic Hausdorff dimension is $1$ by \cite[Corollary 5.1.14]{Adams1996}. The proof of the above theorem is presented in Section \ref{sec:proofs} and it applies to a more general class of Dobi\'nski type sets and all $(a,p)$ Riesz capacities (see Theorem \ref{thm:secondmain}). 
The proof rests on two ideas. One is the use of a discrete/dyadic version of the Riesz capacity. Discrete type capacities have appeared in potential theory in the past (see for example \cite{Benjamini1992, Soardi1994}) and their ``combinatorial'' nature suites very well the dyadic structure of $\cD$. In concrete terms one can show using a {\it recursive formula} (Lemma \ref{lem:recursion}) that  $\cD$ has positive discrete capacity and this, through a comparison theorem for discrete and Riesz capacities, \cite[Theorem 1]{Arcozzi2013} allow one to deduce that $\cD$ has positive Riesz $(\frac{1}{p},p) $-capacity for $1 < p\leq 2$ and zero capacity when $p>2$. Finally we prove a ``Kolmogorov $0-1$'' type lemma (Lemma \ref{lem:01law}) from which we can deduce that in fact $\cD$ is of {\it full} capacity, i.e. $R_{\frac{1}{p},p}(\cD) =  R_{\frac{1}{p},p}(\bT )$ when $1< p \leq 2$. It is worth noticing that the same phenomenon ($0-1$ type law) appears in the study of logarithmic capacity of uniform $G_{\delta}$-sets \cite[Theorem 1.2]{Kleptsyn2021}.

\section{Trees, dyadic capacity and the recursion formula}

Let $T:=\{0,1\}^*$ the free monoid generated by the language $\{0,1\}$ with neutral element $e$. In this context we shall call $T$ the {\it dyadic tree}. The length of a word $x$ is denoted  $|x|$. For two words $x,y\in T$ we denote the largest common prefix of $x$ and $y$ by $x\wedge y$. If $x\wedge y = x$ we write  $x\leq y$. Finally we use the notation $x_-=x0, \, x_+=x1.$ The (Poisson) boundary $\partial T$ of $T$ can be identified with the metric space $\{0,1\}^\bN$ equipped with the metric \[ d(x,y):=2^{-|x\wedge y|}. \]
We will write $\overline{T}:=\partial T \cup T.$ There exists a natural mapping from $\partial T $ to $[0,1]$,
\begin{align}
    \Lambda: \partial T \mapsto [0,1], \quad \Lambda(w_1w_2 \cdots ) = (0.w_1w_2\cdots)_2 
\end{align}
which is onto and Lipschitz continuous. Moreover, every $x\in[0,1]$ which is not a dyadic rational has a unique pre-image. Dyadic rationals have two pre-images under $\Lambda.$  

Our next goal is to develop a potential theory on $\partial T$ which parallels the one we have already seen in $\bT$. A more detailed exposition of the potential theory on the boundary of the tree can be found in \cite{Arcozzi2013, Chalmoukis2019, Arcozzi2021}. Here we shall present only the elements that are essential for our problem.

Let $\varphi$ a non negative function defined on $T$.  The potential of $\varphi$ is given by 
\[ I\varphi (x) := \sum_{y \leq x}\varphi(y), \quad x\in \overline{T}.  \]
  Let $\pi$ be a positive weight function defined on $T$. Then for a set $E \subseteq \partial T$ we define its $\pi$ (discrete) capacity as follows 
\[ \capa_{\pi}(E):=\inf \big\{ \sum_{x\in T} \varphi(x)^p \pi(x) : \varphi \geq 0\, \text{and}\,\,  I\varphi \geq 1 \,\text{on}\, E  \big\}. \]
When $\pi(x)=2^{-|x|(1-ap)}$ we shall refer to the capacity $\capa_\pi=:\capa_{a,p}$ as discrete $(a,p)$ capacity. The relation between the Riesz and discrete capacities can be made explicit. This has been first noted in \cite{Benjamini1992} and generalized further in \cite{Arcozzi2013}.
 
\begin{thm}\cite[Theorem 1]{Arcozzi2013} \label{thm:comparison}Let $p>1, 0<a\leq 1/p$. There exists a constant $c=c(a,p)>0$ such that for any compact set $K\subseteq \partial T$ 
\[ c^{-1}\capa_{a,p}(K) \leq R_{a,p}(\Lambda(K)) \leq c \capa_{a,p}(K). \]
\end{thm}
In fact the restriction that $K$ should be a compact set can be relaxed considerably. By Choquet's capacitability theorem \cite[Theorem 2.3.11]{Adams1996} , Theorem \ref{thm:comparison} holds for all Suslin sets, in particular for all Borel sets.

Discrete capacities satisfy a recursive formula, which is of fundamental importance for our computations. It relates the capacity of a set to the capacities of the parts of its dyadic decomposition. Let $x\in T$ and $E \subseteq \partial T.$ Let also $E_x:=\{w\in \partial T : xw \in E \}$ and $\pi_x(w)=\pi(xw)$. Then we define 
\[ \capa_\pi(E,x)=\capa_{\pi_x}(E_x). \]
Informally, $\capa_\pi(E,x)$ is the capacity of the portion of $E$ that stays below $x$ ``viewed'' from the root $x$.

\begin{thm}\cite[Theorem 30]{Arcozzi2013}\label{lem:recursion}
Let $E\subseteq \partial T$ a Borel set. For every $x\in T$ the following equality holds
\begin{equation}\label{eq:recursion}
    \capa_{\pi}(E,x) = \frac{\capa_{\pi}(E,x_-)+\capa_{\pi}(E,x_+)}{\Bigg(1+\Bigg[\dfrac{ \capa_{\pi}(E,x_-)+\capa_{\pi}(E,x_+)}{\pi(x)} \Bigg]^{p'-1} \Bigg)^{p-1}}. 
\end{equation}
\end{thm}

Finally let us introduce a more general class of Dobi\'nski type sets on the boundary of the tree. This is a rather natural generalization of the set $\cD$. Suppose that $\kappa_n$ is a sequence of positive integers. Let 

\[ D(n,\kappa_n):=\{w\in\partial T : w_{n+1} = w_{n+2} = \dots w_{n+\kappa_n} = 0 \}. \] We define the {\it Dobi\'nski} type set associated to $\kappa_n$ as
\[ D:= \limsup_{n \to \infty} D(n,\kappa_n)  \]
 Notice that if we consider the set $\Lambda(\cup_m D_m)$ where $D_m$ is the Dobi\'nski type set corresponding to the sequence $\kappa_n=2^nm^{-1}$, we obtain ``one half'' of the Dobi\'nski set $\cD$. The other half is obtained by considering the same construction, where instead of ``strings of 0's'' in the definition of $D(n,\kappa_n)$ we consider strings of $1$'s. 
\section{Proof of the main result}\label{sec:proofs}

\begin{lem}\label{lem:01law}
Let $p>1$ and $a>0$ such that $ap\leq 1$ and $E\subseteq \bT$ a Borel set which is invariant under rotations by angles $\theta$, where $\theta$ is a dyadic rational number. Then, either $R_{a,p}(E)= R_{a,p}(\bT) $ or $R_{a,p}(E)= 0$.
\end{lem}

\begin{proof}Assume that $R_{a,p}(E) \neq 0$. By Theorem \cite[Theorem 2.3.10]{Adams1996} there exists a unique non negative function $f_E\in L^p(\bT)$ such that 
\begin{equation}\label{eq:potential}
 \cI_a(f_E) \geq 1, \quad  R_{a,p}\text{-q.e. on } E, \end{equation}
 and 
 \begin{equation}\label{eq:energy}
     \int_{\bT} f_E^p(x)dx = R_{a,p}(E).
 \end{equation}
Let $\theta$ a dyadic rational and define $\rho_\theta f(x) : = f(e^{2\pi i \theta} x )$. Then it is clear that $ \cI_a(\rho_\theta f_E) = \rho_\theta (\cI_a(f_E)) \geq 1 $ on $e^{2\pi i \theta} E = E.$ Therefore $\rho_\theta f_E$ satisfies equations \eqref{eq:potential} and \eqref{eq:energy} and by uniqueness $\rho_\theta f_E = f_E$. Since $\theta$ is dense in $[0,1]$ a calculation with the Fourier coefficients of $f_E$ shows that $f_E = c$ Lebesgue a.e. on $\bT$ for some positive constant $c$.

By equation \eqref{eq:energy} and the fact that $0<R_{a,p}(E)\leq R_{a,p}(\bT)$  we get $0<c^p\leq R_{a,p}(\bT)$. Finally let $y_0\in E$, such that  \eqref{eq:potential} holds;
\[ 1 \leq \cI_a(f_E)(y_0) = c \int_\bT \frac{dx}{|y_0-x|^{1-a}} = c R_{a,p}(\bT)^{-\frac{1}{p}}. \]
\end{proof}

We now turn to the main theorem. The calculation can be carried out for general $(a,p)$, 

\begin{thm}\label{thm:secondmain}
Let $ D $ a Dobi\'nski set associated to a sequence $\kappa_n$ and $a> 0, 1<p<\infty$. In the case  $ap=1$ we have
\begin{itemize}
    \item [(i)]If $\limsup_{n \to \infty}\kappa_n^{-(p-1)}2^{n} > 0$, then $\capa_{\frac{1}{p},p}(D)>0. $
    \item[(ii)]  If $\sum_{n = 1}^\infty \kappa_n^{-(p-1)}2^{n}<  \infty $, then $\capa_{\frac{1}{p},p}(D)=0$.
\end{itemize}
While in the $ap<1$ we have
\begin{itemize}
\item [(a)]If $\limsup_{n\to \infty} (apn-(1-ap)\kappa_n) > -\infty $ then $\capa_{a,p}(D) >0. $,
\item[(b)] $\sum_{n=0}^\infty 2^{  ap n - (1-ap) \kappa_n } < +\infty  $ then $\capa_{a,p}(D)=0$.
\end{itemize}

\end{thm}
\begin{proof} Let $D(n,\kappa_n)$ as before. We start with deriving an exact formula for the discrete  $(a,p)$ capacity of the set $D(n,\kappa_n)$ using the recursive formula (equation \eqref{eq:recursion}). For a positive parameter $r>0$ define the function 
\[ \Phi_r(x):= \frac{x}{(1+rx^{p'-1})^{p-1}}, \,\,\, x>0. \]

An elementary computations shows that the following semigroup law is satisfied 
\[ \Phi_r \circ \Phi_s (x) = \Phi_{r+s}(x), \,\, \forall r,s,x>0. \]
Next we apply $n+\kappa_n$ times the recursive formula \eqref{eq:recursion} for the set $D(n,\kappa_n)$. In the following $c:=\capa_{a,p}(\partial T)$. If we use the symbol $\prod$ for repeated composition of functions  we have \begin{align}\label{eq:composition}
    \capa_{a,p}(D(n,\kappa_n))=2 ^n \prod_{m=1}^n \Phi_{2^{(p'-1)[(n+1-m)+(m-1)(1-ap)]}}\circ \prod_{m=n+1}^{n+\kappa_n} \Phi_{2^{(p'-1)(m-1)(1-ap)}} (c) = 2^n \Phi_{\sigma}(c)
\end{align}
where we have used the fact that $\Phi_r(2x)=2\Phi_{2^{p'-1}r}(x)$ and $\sigma $ is given by
\begin{align*} \sigma = & \sum_{m=1}^n 2^{(p'-1)[(n+1-m)+(m-1)(1-ap)]} + \sum_{m=n+1}^{n+\kappa_n} 2^{(p'-1)(m-1)(1-ap)} \\
 = & \frac{2^{n(p'-1)}-2^{n(1-ap)(p'-1)}}{1-2^{-ap(p'-1)}}+ \frac{2^{(n+\kappa_n)(1-ap)(1-p')}-2^{n(1-ap)(p'-1)}}{2^{(1-ap)(p'-1)}-1}, \,\,\, \text{when} \,\,\, ap<1, \\
 \sigma = &  \frac{2^{n(p'-1)}-1}{1-2^{1-p'}} + \kappa_n,  \,\,\, \text{when} \,\,\, ap=1.
\end{align*}
Consequently if $ap=1$,
\begin{align*}\capa_{\frac{1}{p},p}(D(n,\kappa_n)) & = \dfrac{2^n c}{\Big(1+ \kappa_n c^{p'-1} + \dfrac{2^{n(p'-1)}-1}{1-2^{1-p'}} c^{p'-1} \Big)^{p-1}} \\
& =  \dfrac{ c}{\Big(2^{n(1-p')}+ \kappa_n 2^{n(1-p')} c^{p'-1} + \dfrac{1-2^{n(1-p')}}{1-2^{1-p'}} c^{p'-1} \Big)^{p-1}}.
\end{align*}
From which is easily verified that there exists a constant $A>0$ such that 
\[ \frac{1}{A}\capa_{\frac{1}{p},p}(D(n,\kappa_n)) \leq \kappa_n^{-(p-1)} 2^{n} \leq A \capa_{\frac{1}{p},p}(D(n,\kappa_n)), \quad \forall n\in \bN.  \]
Similarly when $ap<1$,

\begin{align*}\capa_{\frac{1}{p},p}(D(n,\kappa_n)) & = \dfrac{2^n c}{\Bigg(1+\Big[\dfrac{2^{n(p'-1)}-2^{n(1-ap)(p'-1)}}{1-2^{-ap(p'-1)}}+ \dfrac{2^{(n+\kappa_n)(1-ap)(p'-1)}-2^{n(1-ap)(p'-1)}}{2^{(1-ap)(p'-1)}-1}\Big]c^{p'-1} \Bigg)^{p-1}} \\
& =  \dfrac{ c}{\Big(2^{-n(p'-1)}+ \Big[ \dfrac{1-2^{-nap(p'-1)}}{1-2^{-ap(p'-1)}} +\dfrac{2^{-n ap(p'-1)+\kappa_n(1-ap)(p'-1)}-2^{-nap(p'-1)}}{2^{(1-ap)(p'-1)}-1}\Big]c^{p'-1} \Big)^{p-1}}.
\end{align*}
Hence, for some constant $A>0$, 
\[\frac{1}{A} \capa_{\frac{1}{p},p}(D(n,\kappa_n)) \leq  2^{ ap n - (1-ap)\kappa_n} \leq A  \capa_{\frac{1}{p},p}(D(n,\kappa_n)),\quad \forall n \in \bN.\]

The theorem then follows from the estimates which hold for all $ap \leq 1$. We can estimate the capacity of $D$ from above using subadditivity 
\begin{equation*}
    \capa_{a,p}(D)= \capa_{a,p}\Big( \bigcap_{m=1}^\infty \bigcup_{n=m}^\infty D(n,\kappa_n) \Big) \leq \sum_{n=m}^\infty \capa_{a,p}D(n,\kappa_n), \,\,\, \forall m \in \bN.
\end{equation*}
And from below; 
\begin{equation*}
    \capa_{a,p}(D)= \capa_{a,p}\Big( \bigcap_{m=1}^\infty \bigcup_{n=m}^\infty D(n,\kappa_n) \Big) = \lim_m \Big( \capa_{a,p} \bigcup_{n=m}^\infty  D(n,\kappa_n) \Big) \geq \limsup_{n} \capa_{a,p} D(n,\kappa_n).
\end{equation*}
\end{proof}

\begin{proof}[Proof of Theorem \ref{thm:main}]
Let $p>1$. Consider the Dobi\'nski type set $D_m$ corrisponding to the sequence $\kappa_n=2^nm^{-1}$ for some fixed $m$. As noted before $ \Lambda(D_m)\subseteq \cD$. For $1<p\leq2$, by (i) Theorem \ref{thm:main}, we have that $\capa_{\frac{1}{p},p}(D_m)>0$. Therefore the comparison principle,  Theorem \ref{thm:comparison}, gives $ R_{\frac{1}{p},p}(\cD) \geq R_{\frac{1}{p},p}(\Lambda(D_m)) \geq c^{-1} \capa_{\frac{1}{p},p}(D_m) >0.$  Since $\cD$ is invariant under rotations by dyadic rationals, Lemma \ref{lem:01law} implies that $\cD$ has full
capacity.

While when $p>2$, by Theorem \ref{thm:main} (ii), $\capa_{\frac{1}{p},p}(D_m)=0$ for all $m\in \bN.$ By countable subadditivity and the comparison theorem,  $\capa_{\frac{1}{p},p}(\cup_m D_m) = R_{\frac{1}{p},p}(\Lambda(\cup_m D_m)) = 0 $. The ``other half'' of Dobi\'nski's set has zero capacity for the same reason, which concludes the proof. 
\end{proof} 

It is clear that Theorem \ref{thm:main} leaves a gap between the  sufficient conditions for a set to have zero or positive capacity. Although for the Dobi\'nski set this turns out not to be a problem, one can construct Dobi\'nski type sets for which Theorem \ref{thm:main} does not provide a answer about their capacity. We conjecture that the series condition in Theorem \ref{thm:main} in fact characterizes the vanishing of the capacity of Dobi\'nski type sets. 

\subsection*{Acknowledgements} We would like to thank Alberto Dayan for interesting discussions on the problem.

\bibliographystyle{plain}
\bibliography{literature}
\end{document}